\documentclass[11pt]{amsart}
\usepackage{amssymb}
\usepackage{enumitem}
\usepackage{array}
\usepackage{graphicx}
\usepackage{tikz}
\usepackage{hyperref}
\hypersetup{hidelinks}

\newtheorem{theorem}{Theorem}[section]
\newtheorem{lemma}[theorem]{Lemma}
\newtheorem{proposition}[theorem]{Proposition}
\newtheorem{corollary}[theorem]{Corollary}
\newtheorem{conjecture}[theorem]{Conjecture}
\theoremstyle{remark}

\numberwithin{equation}{section}

\newcommand{\nint}[1]{\|#1\|}
\newcommand{\step}[1]{\par\smallskip\noindent\textit{#1}\enspace}

\begin{document}

\title[The Sudler Products Revisited]{The Sudler Products Revisited: Strict monotonicity of the largest absolute coefficient of $(q;q)_m$}

\author{Alexander Berkovich}
\address{University of Florida, Department of Mathematics, Little Hall 32601, Gainesville FL}
\email{alexb@ufl.edu}
\author{Ali K. Uncu}
\address{University of Bath, Department of Computer Science, Claverton Down, Bath}
\email{aku21@bath.ac.uk}


\begin{abstract}
Let $M_m$ be the largest absolute value of a coefficient of the polynomial $(1-q)(1-q^2)\cdots(1-q^m)$. We prove that $M_{m+1}~>~M_m$ for every $m\ge 25$. To that end, we analyze the Sudler products $\prod_{r\le m}~|1~-~e^{ir\theta}|$. 
\end{abstract}

\maketitle

\section{Introduction}


We start by defining $q$-Pochhammer symbols. For an integer $m\ge 0$ or $\infty$ write
\[
(q;q)_m=\prod_{r=1}^{m}(1-q^r)=\sum_{k=0}^{m(m+1)/2}a_{m,k}q^k,
\]
and let
\[
M_m=\max_k|a_{m,k}|
\]
be the largest absolute value of a coefficient. Euler's pentagonal number theorem \[(q;q)_\infty = \sum_{k=-\infty}^\infty (-1)^k q^{\frac{k(3k+1)}{2}} = 1-q-q^2 +q^5+q^7-q^{12}+\dots\] shows that $M_\infty = 1$. However, for finite $m$ this is not the case. In 1964, Sudler \cite{Sudler1,Sudler2} showed that $M_m$ is unbounded and that $\log M_m\sim Km$ for a constant $K>0$, which we call the \textit{Sudler's constant}. Wright \cite{Wright1,Wright2} estimated the location of the maximum of the associated trigonometric product. The authors studied the location of the maximal coefficient further in \cite{BU2}, and explicit conjectures for the location of the absolute maximums are given, with extensive computations up to $m=75000$.

In \cite{BU1}, the authors proved that $(q;q)_m$ has all its coefficients in $\{-1,0,1\}$ exactly when $m\in\{0,1,2,3,5,\infty\}$. More generally, they considered for $h\ge 1$ the sets
\[
\mathcal{S}_h=\{m\ge 0:\ M_m=h\}.
\]
They computed $\mathcal{S}_h$ for $h\le 40$ and stated the following conjecture \cite[Conjecture~4.1]{BU1}.

\begin{conjecture}\label{conj}
For $h>16$, either $\mathcal{S}_h=\emptyset$ or $\mathcal{S}_h=\{i(h)\}$ for a positive integer $i(h)$, and $i(h_1)>i(h_2)$ whenever $h_1>h_2>16$. Moreover, for $h>5$ the set $\mathcal{S}_1\cup \mathcal{S}_2\cup\cdots\cup \mathcal{S}_h$ consists of all integers from $0$ up to some $M(h)$.
\end{conjecture}

In this paper, we prove the following theorem, which implies Conjecture~\ref{conj} as its corollary.

\begin{theorem}\label{thm:main}
$M_{m+1}>M_m$ for every $m\ge 25$.
\end{theorem}

\begin{corollary}\label{cor:conj}
Conjecture~\ref{conj} is true.
\end{corollary}





Our analysis and proofs will rest on the integral representation
\begin{equation}\label{eq:akPm}
a_{m,k}=\frac{1}{2\pi}\int_{-\pi}^{\pi}P_m(\theta)e^{-ik\theta}\,d\theta,\qquad P_m(\theta)=\prod_{r\le m}(1-e^{ir\theta}).    
\end{equation}

It is easy to check the values of $M_m$ for $m\le 2000$ using computer algebra. The real aim of this paper is to show
\[
M_m=S_m\bigl(1+\vartheta_m\bigr),\qquad -0.005533\le\vartheta_m\le 0.004148,
\]
where $S_m = -(e^{\varphi(\theta^*)}/\pi) \sqrt{{2\pi}/{\varphi''(\theta^*)}}$ is an explicit saddle point expression (Section~\ref{sec:local}) with $\varphi(\theta) = \log|P_m(\theta)|$ and $\theta^*$ being the location of $\varphi$'s maximum. We will show that
\begin{equation}\label{eq:Sm_ratio}
S_{m+1}/S_m\ge 1.1 .
\end{equation}
This will yield $M_{m+1}/M_m\ge 1.089$ for all $m\ge 2000$ and imply Theorem~\ref{thm:main} when combined with the exact calculations of $25\leq m\leq2000$.

The modulus $|P_m(\theta)|$ is what we will call \textit{the Sudler product}. Its maximum on the circle sits at $\theta^*\approx 2\pi x_0/m$ with $x_0 \approx 0.79122$, where it is $\approx e^{Km}$ as found in \cite{Sudler2,Wright1}. 

The organization of the paper is as follows. Section~2 is reserved for preliminary results. In Sections~\ref{sec:saddle} and~\ref{sec:local}, we study values near the maximum with an explicit Gaussian approximation of $P_m$ in a window of width $m^{-3/2}$ around $\pm\theta^*$, yielding two-sided bounds for $M_m$. The values away from the maximum will be studied in Section~\ref{sec:away}, an explicit bound $|P_m(\theta)|\le e^{-\gamma m}|P_m(\theta^*)|$ for $2\pi/m\le|\theta|\le\pi$. 
We then prove \eqref{eq:Sm_ratio} in Section~\ref{sec:proof}. The exact value calculations, certification of necessary constants in interval arithmetic, and formalization of the results are discussed in Section~\ref{sec:comp}.

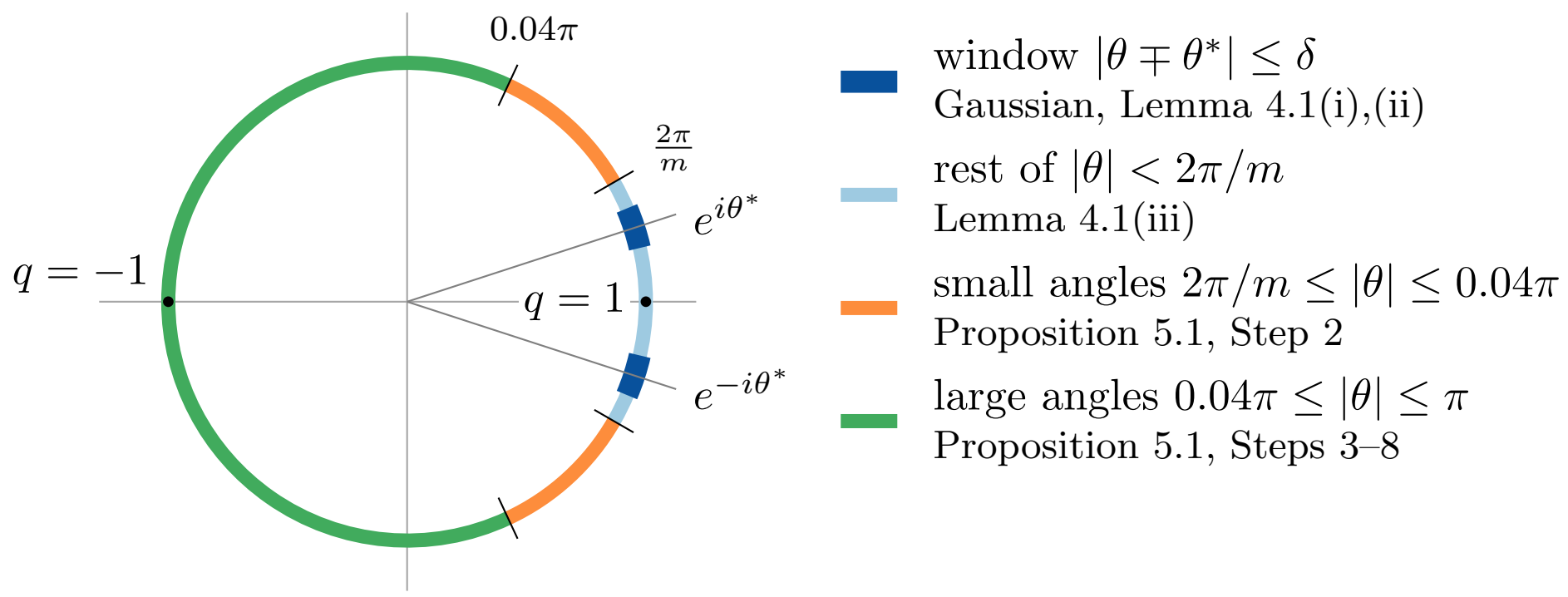
\begin{figure}[tbp]
\centering
\begin{tikzpicture}[scale=1.0]
\definecolor{cwin}{HTML}{08519C}
\definecolor{ccen}{HTML}{9ECAE1}
\definecolor{csma}{HTML}{FD8D3C}
\definecolor{clar}{HTML}{41AB5D}
\def\R{1.9}
\draw[gray!70,thin] (-2.45,0) -- (2.3,0);
\draw[gray!70,thin] (0,-2.3) -- (0,2.3);
\draw[clar,line width=3.2pt] (65:\R) arc (65:180:\R);
\draw[clar,line width=3.2pt] (-65:\R) arc (-65:-180:\R);
\draw[csma,line width=3.2pt] (30:\R) arc (30:65:\R);
\draw[csma,line width=3.2pt] (-30:\R) arc (-30:-65:\R);
\draw[ccen,line width=3.2pt] (-30:\R) arc (-30:30:\R);
\draw[cwin,line width=5pt] (13:\R) arc (13:23:\R);
\draw[cwin,line width=5pt] (-13:\R) arc (-13:-23:\R);
\foreach \a in {30,-30,65,-65} \draw[thin] (\a:\R-0.18) -- (\a:\R+0.18);
\node[left,font=\small,fill=white,inner sep=1pt] at (\R-0.12,0) {$q=1$};
\fill (\R,0) circle (1.2pt);
\fill (-\R,0) circle (1.2pt);
\node[left] at (-\R-0.02,0.22) {\small $q=-1$};
\draw[thin,gray] (0,0) -- (18:\R+0.35) node[right,black] {\small $e^{i\theta^*}$};
\draw[thin,gray] (0,0) -- (-18:\R+0.35) node[right,black] {\small $e^{-i\theta^*}$};
\node at (30:\R+0.55) {\scriptsize $\tfrac{2\pi}{m}$};
\node at (65:\R+0.5) {\scriptsize $0.04\pi$};
\begin{scope}[shift={(3.45,1.75)}]
\def\dy{0.9}
\draw[cwin,line width=5pt] (0,0) -- (0.45,0);
\node[right,align=left,font=\small] at (0.6,0) {window $|\theta\mp\theta^*|\le\delta$\\[-1pt]\footnotesize Gaussian, Lemma~\ref{lem:gauss}(i),(ii)};
\draw[ccen,line width=3.2pt] (0,-\dy) -- (0.45,-\dy);
\node[right,align=left,font=\small] at (0.6,-\dy) {rest of $|\theta|<2\pi/m$\\[-1pt]\footnotesize Lemma~\ref{lem:gauss}(iii)};
\draw[csma,line width=3.2pt] (0,-2*\dy) -- (0.45,-2*\dy);
\node[right,align=left,font=\small] at (0.6,-2*\dy) {small angles $2\pi/m\le|\theta|\le 0.04\pi$\\[-1pt]\footnotesize Proposition~\ref{prop:away}, Step 2};
\draw[clar,line width=3.2pt] (0,-3*\dy) -- (0.45,-3*\dy);
\node[right,align=left,font=\small] at (0.6,-3*\dy) {large angles $0.04\pi\le|\theta|\le\pi$\\[-1pt]\footnotesize Proposition~\ref{prop:away}, Steps 3--8};
\end{scope}
\end{tikzpicture}
\caption{How the circle method is executed. The circle of integration $q=e^{i\theta}$, in the representation $a_{m,k}$ defined by \eqref{eq:akPm}. 
}
\label{fig:circle}
\end{figure}

\begin{figure}[tbp]
\centering
\includegraphics[width=\textwidth]{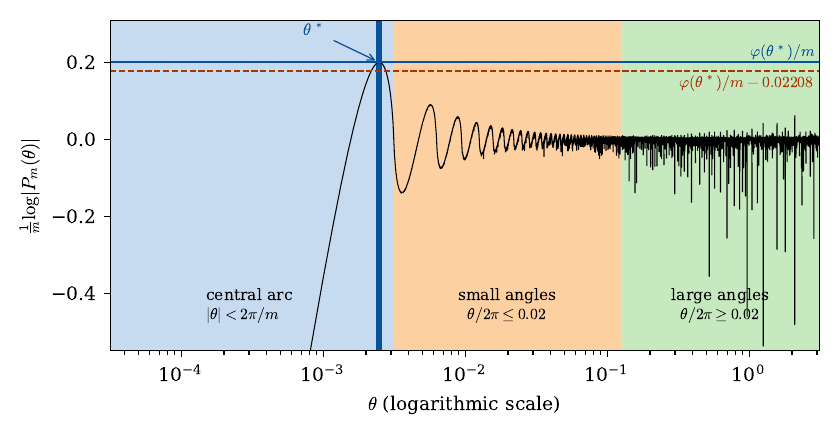}
\caption{The function $\frac1m\log|P_m(\theta)|$ for $m=2000$ and $0<\theta\le\pi$, on a logarithmic $\theta$-axis, with the regions of Figure~\ref{fig:circle}. The dark band is the window $|\theta-\theta^*|\le\delta$. Proposition~\ref{prop:away} asserts that on the small-angle and large-angle arcs the values are bounded above by the dashed line $\varphi(\theta^*)/m-0.02208$. 
The downward spikes come from zeros of $P_m$.
}
\label{fig:split}
\end{figure}

Figure~\ref{fig:circle} shows how the circle of integration $|q|=1$ is divided, and Figure~\ref{fig:split} shows $(1/m)\log|P_m|$ on the different arcs for $m=2000$.

The error bounds that appear in our calculations are non-increasing in $m$. Therefore, it is enough to evaluate it at $m=2000$ to get the necessary bounds that yield the proof of Theorem~\ref{thm:main} when combined. 




\section{Preliminaries}\label{sec:2}





For real $t$ let $\nint{t}$ denote the distance from $t$ to the nearest integer, and put
\[
L(t):=\log|2\sin\pi t|,\qquad\text{and}\qquad \Lambda(x):=\int_0^xL(t)dt,
\]
Note that $L(t)$ is even and $1$-periodic. Its integral over a full period $\int_k^{k+1}L(t) dt$ is 0. $L(t)$ is increasing on $(0,\frac12]$ and decreasing on $[\frac12,1)$. In particular, $L(t)>0$ for $t\in (\frac16,\frac56)$.
Moreover
\begin{equation}\label{eq:sinbounds}
4\nint{t}\le|2\sin\pi t|\le 2\pi\nint{t},\qquad |L'(t)|=\pi|\cot\pi t|\le\frac{1}{\nint{t}}.
\end{equation}
For every integer $k\ge 1$ and real $t$, we have
\begin{equation}\label{eq:product}
\sum_{j=0}^{k-1}L\Bigl(t+\frac{j}{k}\Bigr)=L(kt),
\end{equation} which one can easily get by first writing $L(t) = \log |e^{2\pi i t}-1|$. 
For $\alpha\in\mathbb{R}$, we write
\begin{equation}\label{eq:sigmamalpha}
\Sigma_m(\alpha):=\sum_{r=1}^{m}L(r\alpha)=\log|P_m(2\pi\alpha)|.
\end{equation}
Observe that $\Lambda$ is also $1$-periodic. Since $L(t)\le 0$ on $[0,\frac16]\cup[\frac56,1]$ and $L(t)>0$ in $(\frac16,\frac56)$. Whence,
\[
\max_{x\in\mathbb{R}}\Lambda(x)=\Lambda(\tfrac56).
\]



\begin{lemma}\label{lem:sudler}
For $0\le U<\pi$ and $y>0$, let
\[
F(U)=\int_0^Uu\cot u\,du \qquad \text{and}\qquad f(y)=\frac{\Lambda(y)}{y}.
\]
$F$ has a unique zero $U_0=\pi x_0$ in $(0,\pi)$, and the \textit{Sudler's constant} is defined as follows,
\[
K:=\sup_{y>0}f(y)=f(x_0)=L(x_0).
\]
Furthermore, these constants satisfy the following properties:
\[
x_0\in[0.791224,0.791229],\quad K\in[0.198607,0.198629],\quad \Lambda(\tfrac56)\le 0.161538.
\]
\end{lemma}

\begin{proof}
\step{The range $0<y<1$.} Integration by parts gives $\Lambda(y)=yL(y)-\frac1\pi F(\pi y)$, hence
\[
f'(y)=\frac{F(\pi y)}{\pi y^2}.
\]
Since $F'(U)=U\cot U$, the function $F(U)$ is positive and increasing on $(0,\frac\pi2]$ and decreasing on $[\frac\pi2,\pi)$. As $U\to\pi$ we have $F(U)\to-\infty$. So $F(U)$ has a unique zero $U_0\in(\frac\pi2,\pi)$, $f(y)$ attains its maximum over $(0,1)$ at $x_0=U_0/\pi$, and there $\Lambda(x_0)=x_0L(x_0)$.

\step{The range $y\ge 1$.} Here $f(y)=\Lambda(y-\lfloor y\rfloor)/y\le \Lambda(\tfrac56)$, where $\lfloor\cdot\rfloor$ is the classical floor function.

We leave the comments on numerical computations to Section~\ref{sec:comp}. In particular $\Lambda(\tfrac56)<K$, so the supremum over $y>0$ is $f(x_0)$.
\end{proof}

Now we present four results that will be used in bounding the errors in later computations.

\begin{lemma}\label{lem:riemann}
Let $g$ have bounded variation $V_0^U(g)$ on $[0,U]$, let $m\ge 1$ and $\Delta=U/m$. Then
\[
\Bigl|\sum_{r=1}^{m}g(r\Delta)-\frac1\Delta\int_0^Ug(u)\,du\Bigr|\le V_0^U(g).
\]
\end{lemma}

\begin{proof}
Clearly, the difference equals
\[
\frac1\Delta\sum_r\int_{(r-1)\Delta}^{r\Delta}\bigl(g(r\Delta)-g(u)\bigr)\,du,
\]
and each integrand is bounded by the variation of $g$ on $[(r-1)\Delta,r\Delta]$.
\end{proof}

The following bounds $L(t)$ will be used repeatedly.

\begin{lemma}\label{lem:onesided}
Let $0<\Delta\le\frac1{50}$, $n\ge 1$, 
and $t_r=r\Delta \not\in \mathbb{Z}$.
\begin{enumerate}[label=\textup{(\alph*)}]
\item We have
\[
\sum_{r=1}^{n}L(t_r)\le\frac1\Delta\Lambda(n\Delta)+(n\Delta+2)\Bigl(\log\frac{1}{2\Delta}+1\Bigr)+(2y+3)\log 2.
\]
\item If $R\subseteq\{1,\dots,n\}$ contains every $r$ with $\nint{t_r}\ge\Delta$, then
\[
\sum_{r\in R}L(t_r)\le\frac1\Delta\Lambda(n\Delta)+3(n\Delta+2)\Bigl(\log\frac{1}{6\Delta}+1\Bigr)+(2y+3)\log 2.
\]
\item If $n\ge 19$ and $\Delta=x_0/n$, then
\[
\sum_{r=1}^{n}L(t_r)\ge nK-\log 2.
\]
\end{enumerate}
\end{lemma}

The left sides sums diverge to $-\infty$ if $t_r$ is an integer. This makes the first two cases trivial and the third case becomes false. Hence, we exclude this case in the statement.

\begin{proof}
We call $r$ rising if the fractional part $\{t_r\}$ lies in $(0,\frac12]$, falling if it lies in $(\frac12,1)$, and \emph{exceptional} if $|\{t_r\}-\frac12|<\Delta$.

For a rising non-exceptional $r$ let $C_r=[t_r,t_r+\Delta]$, and for a falling non-exceptional $r$ let $C_r=[t_r-\Delta,t_r]$. In both cases $C_r$ lies in an interval on which $L$ is monotone, and
\[
L(t_r)\le\frac1\Delta\int_{C_r}L(t) dt,
\]
because $L(t_r)$ is the minimum of $L(t)$ on $C_r$.

It is easy to see that the $C_r$'s of non-exceptional $r$ have disjoint interiors. Indeed, $C_r$ and $C_{r+1}$ could overlap only if $r$ is rising and $r+1$ is falling. Then $t_r\le k+\frac12<t_{r+1}$ for some integer $k$, which makes $r$ exceptional.

If there is a point in $[0,n\Delta] \backslash\bigcup_{r\in R}C_r$, then it lies in one of two kinds of gap.
First type, a gap $[t_r,t_{r+1}]$ with an exceptional end. Such a gap lies within $2\Delta$ of a half-integer, where $L(t)>0$, because $2\Delta<\frac13$.
Or, the other type, where we have a gap around an integer $k$ (for $k=0$ this is $[0,t_1]$). In case (a) it is a single interval $[t_r,t_{r+1}]$ of length $\Delta$ containing $k$. In case (b), together with the cells of the dropped points $r\notin R$, which satisfy $\nint{t_r}<\Delta$, it lies in an interval $J_k$ of length $3\Delta$ containing $k$.

There are at most $n\Delta+2$ integers $k$ involved in the calculations. Since $L(t)<0$ on $J_k\subset(k-\frac16,k+\frac16)$, we use \eqref{eq:sinbounds} together with the fact that $\int_J\log(4|t|)\,dt$, over intervals of fixed length $\ell$ containing $0$, is smallest for the centred one. This gives
\[
-\frac1\Delta\int_{J_k}L\le-\frac1\Delta\int_{-\ell/2}^{\ell/2}\log(4|t|)\,dt=\frac{\ell}{\Delta}\Bigl(\log\frac{1}{2\ell}+1\Bigr),
\]
with $\ell=\Delta$ in case (a) and $\ell=3\Delta$ in case (b).


If the cells stick out of $[0,n\Delta]$, that can be at most by $\Delta$, which contributes at most $\log 2$. There are at most $2(n\Delta+1)$ exceptional $r$, each contributing $L(t_r)\le\log 2$. Adding up gives (a) and (b).

{Part (c):} Compare a rising $t_r$ with the cell $[t_r-\Delta,t_r]$ and a falling one with $[t_r,t_r+\Delta]$, so that now $L(t_r)\ge\frac1\Delta\int L$ over the cell. Since $t_n=x_0$ and $x_0+\Delta<\frac56$, these cells are disjoint and cover $[0,x_0+\Delta]$, except for one cell straddling $\frac12$ on which $L\le\log 2$. Since $L>0$ on $[x_0,x_0+\Delta]$, we get
\[
\sum_rL(t_r)\ge\frac1\Delta\Lambda(x_0)-\log 2=nK-\log 2.\qedhere
\]
\end{proof}


The next two lemmas bound sums of $L$ over $r$ points spread evenly around the circle. They are the main tools for large angles in Section~\ref{sec:away}.

\begin{lemma}\label{lem:grid}
Let $r\ge 3$, $a\in\mathbb{R}$, and let $y_j\in[a+\frac jr,a+\frac{j+1}r)$ for $0\le j<r$. Then
\[
\sum_{j=0}^{r-1}L(y_j)\le\log r+\log\frac{\pi^2}{2}.
\]
\end{lemma}

\begin{proof}
We work modulo $1$. The grid $G=a+\frac1r\mathbb{Z}$ cuts the interval (which can also be viewed as a circle thanks fo the modulo 1 computations) into $r$ cells, one containing each $y_j$. Let $J_0$ be the cell containing $0$ and $J_h$ the cell containing $\frac12$. Write $J_0=[-t,s)$ with $s,t\ge 0$ and $s+t=\frac1r$. By continuity in $a$ we may assume $s,t>0$. Let $g_1=s$ and $g_r=-t$ be the grid points adjacent to $0$.

The cells strictly between $J_0$ and $J_h$ on the right of $0$ lie in $(0,\frac12]$, where $L$ increases, so their points satisfy $L(y_j)\le ~L(\text{right endpoint})$. The cells strictly between $J_h$ and $J_0$ lie in $[\frac12,1)$, where $L$ decreases, so $L(y_j)\le L(\text{left endpoint})$. In this way every grid point except $g_1$ and $g_r$ is used exactly once.

We use $L(y_j)\le\log 2$ on $J_h$ and, by \eqref{eq:sinbounds}, $L(y_j)\le\log(2\pi\max(s,t))$ on $J_0$. By \eqref{eq:product} and \eqref{eq:sinbounds},
\[
\sum_{g\in G}L(g)=L(rs)\le\log\bigl(2\pi r\min(s,t)\bigr),\qquad L(g_1)+L(g_r)\ge\log(16st).
\]

Hence, altogether
\begin{align*}
\sum_jL(y_j)&\le\log\bigl(2\pi r\min(s,t)\bigr)-\log(16st)+\log\bigl(4\pi\max(s,t)\bigr)=\log\frac{\pi^2r}{2}.\qedhere
\end{align*}
\end{proof}

\begin{lemma}\label{lem:perturb}
Let $r\ge 2$, $0\le\eta\le\frac1{4r}$, $x\in\mathbb{R}$, and let $x_0,\dots,x_{r-1}$ satisfy $|x_j-x|\le\eta$. Then
\[
\sum_{j=0}^{r-1}L\Bigl(\frac jr+x_j\Bigr)\le\Lambda_{r\eta}(rx)+E_r(\eta),
\]
where
\[
\Lambda_c(u)=\log\Bigl(|2\sin\pi u|\Bigl(1+\frac{c}{\nint{u}}\Bigr)\Bigr)\quad\text{(extended by $\log(2\pi c)$ at $u\in\mathbb{Z}$)}
\]
and
\[
E_r(\eta)=\eta r\,e_r,\qquad e_r=\sum_{k=1}^{\lceil r/2\rceil}\Bigl(\frac{1}{k-3/4}+\frac{1}{k-1/4}\Bigr).
\]
\end{lemma}

\begin{proof}
Let $z_j=\frac jr+x$ and $d=\nint{rx}/r\le\frac1{2r}$, and let $j_0$ be the index with $\nint{z_{j_0}}=d$.

\step{The points $j\ne j_0$.} The distances $\nint{z_j}$ for $j\ne j_0$ are distinct values of the form
\[
\frac kr-d\ge\frac{k-1/2}{r}\qquad\text{or}\qquad\frac kr+d\ge\frac kr,\qquad 1\le k\le\lceil r/2\rceil.
\]
By \eqref{eq:sinbounds} and the mean value theorem,
\[
L\Bigl(\frac jr+x_j\Bigr)\le L(z_j)+\frac{\eta}{\nint{z_j}-\eta}.
\]
Here $\nint{z_j}-\eta\ge\frac{k-3/4}{r}$ or $\ge\frac{k-1/4}{r}$ respectively. Summing over $j\ne j_0$ gives at most $E_r(\eta)$.

\step{The point $j_0$.} Here $\nint{\frac{j_0}r+x_{j_0}}\le d+\eta\le\frac12$. Using $\sin(a+b)\le\sin a+b\cos a$ and $\pi\cot\pi d\le 1/d$,
\[
|2\sin\pi(d+\eta)|\le|2\sin\pi d|\,(1+\eta/d).
\]
So $L(\frac{j_0}r+x_{j_0})\le L(z_{j_0})+\log(1+\eta/d)$.

\step{Conclusion.} Finally $\sum_jL(z_j)=L(rx)$ by \eqref{eq:product}, and $\eta/d=r\eta/\nint{rx}$. The case $d=0$ follows by continuity.
\end{proof}

\section{The Saddle Point}\label{sec:saddle}

Now we locate the maximum $\theta^*$ of $|P_m|$ on $(0,2\pi/m)$ and bound the first four derivatives of $\log|P_m|$ at that point. 

Fix a positive integer $m$ and write
\[
\varphi(\theta)=\Sigma_m(\theta/2\pi)=\sum_{r=1}^{m}\log|2\sin(r\theta/2)|\qquad(0<\theta<2\pi/m).
\]

\begin{lemma}\label{lem:polar}
For $0<\theta<2\pi/m$ we have
\[
P_m(\theta)=e^{\varphi(\theta)}e^{i(N\theta-m\pi/2)},\qquad N=m(m+1)/4.
\]
Moreover $\varphi''<0$ on $(0,2\pi/m)$, so $\varphi$ has a unique maximum $\theta^*=\theta^*_m$ there.
\end{lemma}

\begin{proof}
For $0<\beta<2\pi$ we have $1-e^{i\beta}=2\sin(\beta/2)\,e^{i(\beta-\pi)/2}$ with $\sin(\beta/2)>0$. Also
\[
\varphi''(\theta)=-\sum_r\frac{r^2}{4\sin^2(r\theta/2)}<0,
\]
and $\varphi\to-\infty$ at both ends of the interval.
\end{proof}


Now we want to bound the derivatives using bounded valuation results we provided in Section~\ref{sec:2} for large $m$.

Let $s(u)=u/\sin u$. Differentiating $\log|2\sin(r\theta/2)|$ and writing $(r/2)^j=(u_r/\theta)^j$ with $u_r=r\theta/2$, we get
\begin{equation}\label{eq:derivs}
\varphi^{(j)}(\theta)=\theta^{-j}\sum_{r=1}^{m}w_j(u_r),\qquad w_j(u)=u^j\frac{d^j}{du^j}\log\sin u,
\end{equation}
where 
\[
w_1=s\cos u,\quad w_2=-s^2,\quad w_3=2s^3\cos u,\quad w_4=-2s^4(1+2\cos^2u).
\]
Note that $w_2,w_4<0$.

\begin{proposition}\label{prop:saddle}
Let $m\ge 2000$, fix $c^*=m\theta^*/2\pi$. 
Then
\begin{enumerate}[label=\textup{(\roman*)}]
\item $0.790624\le c^*\le 0.791829$;
\item $0.129501\,m^3\le -\varphi''(\theta^*)\le 0.131906\,m^3$;
\item $|\varphi^{(3)}(\theta^*)|\le 0.057737\,m^4$;
\item $|\varphi^{(4)}(\theta)|\le 0.852654\,m^5$ whenever $|\theta-\theta^*|\le 4/|\varphi''(\theta^*)|^{1/2}$.
\end{enumerate}
\end{proposition}

\begin{proof}
The basic estimate comes from applying Lemma~\ref{lem:riemann} to \eqref{eq:derivs} with $\Delta=\theta/2$ and $U=m\theta/2$:
\begin{equation}\label{eq:riemannderiv}
\Bigl|\theta^j\varphi^{(j)}(\theta)-\frac2\theta\int_0^Uw_j\Bigr|\le V_0^U(w_j).
\end{equation}

By concavity, for (i), it suffices to show
\[
\varphi'>0\ \text{at}\ U_a=\pi\cdot 0.790624\qquad\text{and}\qquad\varphi'<0\ \text{at}\ U_b=\pi\cdot 0.791829,
\]
where $U=m\theta/2$. For $j=1$ we have $\int_0^Uw_1=F(U)$, and \eqref{eq:riemannderiv} gives
\[
\theta\varphi'(\theta)\ge\frac mUF(U)-V_0^U(w_1).
\]
Interval arithmetic gives $F(U_a)>0$ and $\frac{2000}{U_a}F(U_a)>V_0^{U_a}(w_1)$. The left side increases with $m$, so $\varphi'>0$ at $U_a$ for all $m\ge 2000$. The argument at $U_b$ is symmetric.

With $\theta=2U/m$, \eqref{eq:riemannderiv} becomes, for $j=2$,
\[
-\frac{\varphi''(\theta^*)}{m^3}\in\frac{1}{4U^3}\int_0^U(-w_2)\pm\frac{V_0^U(w_2)}{4U^2m}.
\]
Similarly
\[
\frac{|\varphi^{(3)}|}{m^4}\le\frac{1}{8U^4}\Bigl|\int_0^Uw_3\Bigr|+\frac{V_0^U(w_3)}{8U^3m},\qquad
\frac{|\varphi^{(4)}|}{m^5}\le\frac{1}{16U^5}\int_0^U(-w_4)+\frac{V_0^U(w_4)}{16U^4m}.
\]
Once we replace $1/m$ by $1/2000$, we see that for (ii) and (iii), part (i) allows us to let $U$ range over the interval $\pi\cdot[0.790624,0.791829]$. Similarly, for (iv), the window $|\theta-\theta^*|\le 4/|\varphi''(\theta^*)|^{1/2}$ corresponds to
\[
|U-U^*|\le 2/(0.129501\,m)^{1/2}\le 0.124274,
\] where $U^* =m\theta^*/2$,
so let $U$ range over that widened interval.
\end{proof}

\section{Major Arc Analysis}\label{sec:local}

In this section we approximate $P_m$ by a Gaussian in a small window around $\theta^*$ and deduce bounds for $M_m$ in terms of the saddle point term $S_m$ and the remainder $R_m$.

We first give a list of constants that will be used in the bounds.
Let $\sigma=1/|\varphi''(\theta^*)|^{1/2}$
and put
\begin{gather*}
\delta=4\sigma,\qquad a_3=\frac{|\varphi^{(3)}(\theta^*)|}{6},\qquad b_4=\frac1{24}\max_{|t|\le\delta}|\varphi^{(4)}(\theta^*+t)|,\\
\rho=a_3\delta^3+b_4\delta^4,\qquad\text{and}\qquad \kappa=8-\rho.
\end{gather*}
Define
\begin{gather*}
T_1=3b_4\sigma^4,\qquad T_2=\frac{e^\rho}{2}\bigl(15a_3^2\sigma^6+96\sqrt{\frac{2}{\pi}}a_3b_4\sigma^7+105b_4^2\sigma^8\bigr),\\
T_3=\frac{e^{-8}}{2\sqrt{2\pi}},\qquad\text{and}\qquad T_4=\frac{8 e^{-\kappa}}{\kappa\sqrt{2\pi}}.
\end{gather*}

\begin{lemma}
\label{lem:gauss}
Assume $\theta^*+\delta<2\pi/m$ and $\kappa>0$. For real $j$ let
\[
I(j)=\int_{-\delta}^{\delta}e^{\varphi(\theta^*+t)-\varphi(\theta^*)}e^{ijt}\,dt.
\]
Then
\begin{enumerate}[label=\textup{(\roman*)}]
\item $|I(j)|\le\sqrt{2\pi}\,\sigma\,(1+T_1+T_2)$;
\item $|I(j)-\sqrt{2\pi}\,\sigma e^{-j^2\sigma^2/2}|\le\sqrt{2\pi}\,\sigma\,(T_1+T_2+T_3+3a_3|j|\sigma^4)$;
\item the integral of $e^{\varphi(\theta)-\varphi(\theta^*)}$ over $\{\theta\in(0,2\pi/m):|\theta-\theta^*|>\delta\}$ is at most $\sqrt{2\pi}\,\sigma\,T_4$.
\end{enumerate}
\end{lemma}

\begin{figure}[tbp]
\centering
\includegraphics[width=\textwidth]{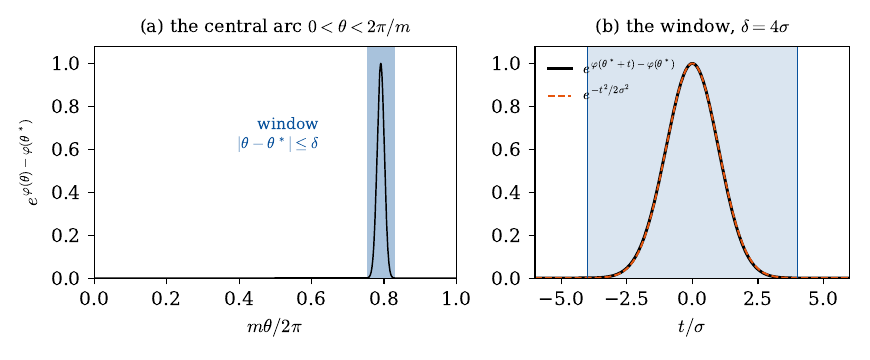}
\caption{The central arc for $m=2000$. (a) The function $e^{\varphi(\theta)-\varphi(\theta^*)}$ on $0<\theta<2\pi/m$; the shaded band is the window $|\theta-\theta^*|\le\delta$, about $8\%$ of the arc ($4\%$ on each side). (b) The window in the variable $t=\theta-\theta^*$, compared with the Gaussian $e^{-t^2/2\sigma^2}$ of Lemma~\ref{lem:gauss}; on the window the two differ by at most $0.0055$. For $m=2000$ one has $c^*=m\theta^*/2\pi=0.79097$ and $1/(\varphi''(\theta^*)m^3)=-0.13076$, inside the ranges of Proposition~\ref{prop:saddle}.}
\label{fig:window}
\end{figure}

Figure~\ref{fig:window} shows the central arc and the window for $m=2000$.

\begin{proof}
Since $\varphi'(\theta^*)=0$, Taylor series of $\varphi(\theta)$ gives
\[
\varphi(\theta^*+t)-\varphi(\theta^*)=-\frac{t^2}{2\sigma^2}+R_3(t)+R_4(t),
\]
where
\[
R_3(t)=\frac{\varphi^{(3)}(\theta^*)t^3}{6},\qquad |R_4(t)|\le b_4t^4.
\]
For $|R_3+R_4|\le\rho$,
\begin{equation}\label{eq:1}
e^{R_3+R_4}=1+R_3+R_4+E,\qquad |E|\le\frac{e^\rho}{2}(a_3|t|^3+b_4t^4)^2.
\end{equation}

Integrate each term of \eqref{eq:1} against $g(t)=e^{-t^2/2\sigma^2}$, using \[\int g\,t^{2k}~=~\sqrt{2\pi}\sigma^{2k+1}(2k-1)!!:\]
The term $1$ gives $\sqrt{2\pi}\sigma e^{-j^2\sigma^2/2}$, minus the mass of $g$ outside the window, which is $\sqrt{2\pi}\sigma\operatorname{erfc}(2\sqrt2)\le\sqrt{2\pi}\sigma T_3$, where $\operatorname{erfc}$ is the classical \textit{complementary error function}.
Since $R_3$ is odd, only $iR_3\sin(jt)$ survives against $e^{ijt}$, and $|\sin jt|\le|j||t|$ gives the term $3a_3|j|\sigma^4$. Without the phase factor this term integrates to $0$, which gives (i).
Together, $R_4$ and $E$ give $T_1$ and $T_2$.

By concavity (Lemma~\ref{lem:polar}), $\varphi$ lies below its chords from $\theta^*$ through $\theta^*\pm\delta$, where $\varphi-\varphi(\theta^*)\le-8+\rho=-\kappa$. Hence the integrand is at most $e^{-\kappa|t|/\delta}$ for $|t|>\delta$, and its integral is at most $2\delta e^{-\kappa}/\kappa$.
\end{proof}

Now we move on to bounding $M_m$.
Define
\[
S_m=\frac{e^{\varphi(\theta^*)}}{\pi}\sqrt{-\frac{2\pi}{\varphi''(\theta^*)}},\qquad R_m=\frac1\pi\int_{2\pi/m}^{\pi}e^{\varphi(\theta)}\,d\theta,
\]
where in $R_m$ we write $\varphi(\theta)=\log|P_m(\theta)|$ for all $\theta$.

\begin{lemma}\label{lem:twosided}
Under the hypotheses of Lemma~\ref{lem:gauss}, put $j^*=\pi/\theta^*+\frac12$. Then
\begin{align*}
M_m&\ge S_m\Bigl(1-\frac{(j^*)^2\sigma^2}{2}-\frac{(\theta^*)^2}{8}-T_1-T_2-T_3-3a_3j^*\sigma^4-T_4\Bigr)-R_m,\intertext{and}
M_m&\le S_m(1+T_1+T_2+T_4)+R_m.
\end{align*}
\end{lemma}

\begin{proof}
Split the integral for $a_{m,k}$ \eqref{eq:akPm} into $|\theta|<2\pi/m$ and the rest. The rest is at most $R_m$ in absolute value.

On $0<|\theta|<2\pi/m$, Lemma~\ref{lem:polar} gives
\[
P_m(\theta)e^{-ik\theta}=e^{\varphi(\theta)}e^{i(j\theta-m\pi/2)},\qquad j=N-k.
\]
The values at $\pm\theta$ are conjugate, so this part equals $\frac1\pi\operatorname{Re}\int_0^{2\pi/m}$. 
Writing $\theta=\theta^*+t$ in the window, it becomes
\[
\frac{e^{\varphi(\theta^*)}}{\pi}\operatorname{Re}\bigl(e^{i\beta}I(j)\bigr)+\hat{O}\bigl(S_mT_4\bigr),\qquad\beta=j\theta^*-\frac{m\pi}{2},
\]
where $\hat{O}(X)$ denotes a quantity of absolute value at most $X$, and the error term comes from Lemma~\ref{lem:gauss}(iii).

Note that $\sqrt{2\pi}\sigma e^{\varphi(\theta^*)}/\pi=S_m$, so Lemma~\ref{lem:gauss}(i) gives the upper bound for every $k$.
As $k$ runs through $\mathbb{Z}$, $j$ runs through $N-\mathbb{Z}$, so $j\theta^*$ runs through an arithmetic progression with difference $\theta^*$. Choose an integer $\ell$ with $|m\pi/2+2\pi\ell|\le\pi$. The closed interval of length $\theta^*$ centred at $m\pi/2+2\pi\ell$ contains some $j\theta^*$. For this $j$ we have
\[
\cos\beta\ge\cos(\theta^*/2)\qquad\text{and}\qquad |j|\le j^*.
\]
Also $k=N-j\in[0,m(m+1)/2]$ because $j^*<N$. Lemma~\ref{lem:gauss}(ii) now gives
\[
\operatorname{Re}(e^{i\beta}I(j))\ge\sqrt{2\pi}\sigma\bigl(\cos(\theta^*/2)e^{-j^2\sigma^2/2}-T_1-T_2-T_3-3a_3|j|\sigma^4\bigr).
\]
Finally use $\cos(\theta^*/2)e^{-x}\ge 1-(\theta^*)^2/8-x$.
\end{proof}

\begin{corollary}\label{cor:twosided}
For every $m\ge 2000$,
\[
S_m(1-0.005533)-R_m\le M_m\le S_m(1+0.004148)+R_m.
\]
\end{corollary}

\begin{proof}
Insert the bounds of Proposition~\ref{prop:saddle} into constants of Section~\ref{sec:local}:
\begin{gather*}
\sigma^2\le(0.129501\,m^3)^{-1},\qquad a_3\le 0.057737\,m^4/6,\qquad b_4\le 0.852654\,m^5/24,\\
2\pi\cdot 0.790624/m\le\theta^*\le 2\pi\cdot 0.791829/m.
\end{gather*}

With these values $\rho$ is a sum becomes a constant multiples of $m^{-1/2}$ and $m^{-1}$, and $\rho\le 0.5667$ at $m=2000$. The numbers $T_1$, $3a_3j^*\sigma^4$, $(j^*)^2\sigma^2/2$ and $(\theta^*)^2/8$ are constant multiples of non-positive powers of $m$ (with $j^*\le m/(2\cdot 0.790624)+1/2$). Not that $T_2$ and $T_4$ decrease as $\rho$ decreases, and $T_3$ does not depend on $m$.

Therefore, every term is non-increasing in $m$, and the values at $m=2000$, computed in interval arithmetic, give the stated constants of the Corollary.

Finally, the hypothesis $\theta^*+\delta<2\pi/m$ holds because $U^*+0.124274<\pi$.
\end{proof}

\section{Away from the saddle}\label{sec:away}

We now show that $|P_m(\theta)|$ is exponentially smaller than its maximum outside the window $|\theta|<2\pi/m$. Small angles are handled by Lemma~\ref{lem:onesided}(a), large angles by a rational approximation and Lemmas~\ref{lem:grid} and~\ref{lem:perturb}.

\begin{proposition}\label{prop:away}
For $m\ge 2000$ and $2\pi/m\le|\theta|\le\pi$,
\[
\log|P_m(\theta)|\le\varphi(\theta^*)-0.02208\,m.
\]
Consequently
\[
R_m/S_m\le\pi\sqrt{-\varphi''(\theta^*)/2\pi}\,e^{-0.02208m}\le 3\cdot 10^{-15}.
\]
\end{proposition}

\begin{proof}
By symmetry we may take $\theta=2\pi\alpha$ with $1/m\le\alpha\le 1/2$. By Lemma~\ref{lem:onesided}(c) and the maximality of $\theta^*$,
\[
\varphi(\theta^*)\ge\Sigma_m(x_0/m)\ge Km-\log 2.
\]
We show
\[
\Sigma_m(\alpha)\le(K-0.02208)m-3\log 2;
\]
this is stronger than needed, since the second term is $-3\log 2$ rather than $-\log 2$.

Apply Lemma~\ref{lem:onesided}(a) with $\Delta=\alpha$ with small angles $\alpha \le 0.02$, $n=m$ and $y=\alpha m\ge\max(1,2000\alpha)$. We use $\Lambda(y)\le \Lambda(\tfrac56)$ and
\[
2/m\le 2\alpha,\qquad 2/m\le 0.001,\qquad 3/m\le 3\alpha,\qquad 3/m\le 0.0015.
\]
This gives $\Sigma_m(\alpha)\le m\,b(\alpha)$ with
\begin{align*}
b(\alpha)&=\frac{\Lambda(\tfrac56)}{\max(1,2000\alpha)}+\min(3\alpha,\alpha+0.001)\Bigl(\log\frac{1}{2\alpha}+1\Bigr)\\
&\qquad+\min(5\alpha,2\alpha+0.0015)\log 2.
\end{align*}
This bound does not depend on $m$. Interval arithmetic over $0<\alpha\le 0.02$ (on $(0,10^{-8}]$ every term is monotone) gives
\[
K-3\log 2/2000-\sup b\ge 0.02208.
\]

For large angles $0.02\le\alpha\le\frac12$ 
by Dirichlet's theorem, there are coprime $p$ and $q$ with
\[
1\le q\le 400,\qquad \beta:=|q\alpha-p|\le 1/401.
\]
Since $\nint{\alpha}\ge 0.02>1/401$, we have $q\ge 2$. As $L$ is even we may replace $\alpha$ by $-\alpha$, so assume $q\alpha-p=\beta\ge 0$, and put $\varepsilon=\beta/q$.

Write $r=bq+s$ with $1\le s\le q$ and $0\le b<B=\lfloor m/q\rfloor$; the remaining fewer than $q$ terms contribute at most $(q-1)\log 2$. In block $b$ the points $r\alpha$ are congruent to $j_r/q+x_r$, where $x_r=r\varepsilon$ and $j_r\equiv rp\pmod q$ runs over all residues. Moreover
\[
|x_r-x_{r'}|\le\eta=(q-1)\varepsilon<\beta.
\]

Case 1: When $60\le q\le 400$ each full block satisfies the hypothesis of Lemma~\ref{lem:grid}, with $a=\min_rx_r$, because $\eta<1/q$. Hence,
\[
\Sigma_m(\alpha)\le\frac mq\Bigl(\log q+\log\frac{\pi^2}{2}\Bigr)+(q-1)\log 2.
\]
Observe that the bound above divided by $m$ is non-increasing in $m$.

Case 2: When $2\le q<60$, we get $\eta<1/401\le 1/(4q)$. Lemma~\ref{lem:perturb}, applied with $x=x_r$ for each $r$ in the block and averaged over the block, gives
\[
\sum_{r\in\text{block}}L(r\alpha)\le\frac1q\sum_{r\in\text{block}}\Lambda_c(r\beta)+E_q(\eta),\qquad c=q\eta=(q-1)\beta,
\]
since $qx_r=r\beta$. We sum over blocks, with $n=Bq\le m$ and $y=m\beta$, and treat the terms as follows.

Near points with $\nint{r\beta}<\beta$ (at most $2(y+2)$ of them) satisfy
\[
\Lambda_c(r\beta)\le\log(2\pi(\beta+c))=\log(2\pi q\beta)<0,
\]
by \eqref{eq:sinbounds}, since $2\pi q\beta\le 2\pi\cdot 59/401<1$.
For the remaining points, $\Lambda_c\le L(r\beta)+c/\nint{r\beta}$. Here
\[
\sum c/\nint{r\beta}\le 2(q-1)(y+2)\bigl(1+\log(\tfrac1{2\beta})\bigr),
\]
since around each of at most $y+2$ integers these points are spaced $\beta$ apart at distances at least $\beta$.
The sum of $L$, $\sum L(r\beta)$ over the remaining points is bounded by Lemma~\ref{lem:onesided}(b) with $\Delta=\beta$, giving $\frac1\beta\Lambda(n\beta)=nf(n\beta)$ plus explicit errors.
In general $nf(n\beta)\le mK$. If $y<\frac1{2\pi}$, then $nf(n\beta)\le(m-q)\log(2\pi y)$, by $f(u)\le L(u)\le\log(2\pi u)$ for $u<\frac12$ and $n\ge m-q$. If $\frac1{2\pi}\le y<\frac12$, then $nf(n\beta)\le m\min(K,\log(2\pi y))$.

Altogether
\begin{multline}\label{eq:caseS}
\Sigma_m(\alpha)\le\frac1q\Bigl[\mathrm{main}+3(y+2)\Bigl(\log\frac1{6\beta}+1\Bigr)+(2y+3)\log 2\\
+2(q-1)(y+2)\Bigl(1+\log\frac1{2\beta}\Bigr)\Bigr]
+\frac mqE_q\Bigl(\frac{q-1}{q}\beta\Bigr)+(q-1)\log 2.
\end{multline}

If $y\le 0.15$, then $\nint{r\beta}\le r\beta\le y$ for every $r\le n$ and $c\le 0.005$, so $\Lambda_c(r\beta)\le\log(2\pi(y+c))<0$ by \eqref{eq:sinbounds}. Hence
\[
\Sigma_m(\alpha)\le\frac mqE_q\Bigl(\frac{q-1}{q}\beta\Bigr)+(q-1)\log 2\le 0.15\,e_q+58\log 2,
\]
which is below the claimed bound.
Otherwise divide \eqref{eq:caseS} by $m$, and bound the main term by $mK$.

For $0.15\le y\le0.5$, at any fixed $y$ each term is non-increasing in $m\ge 2000$: the terms are of the form $(\log m+O(1))/m$ or $1/m$; $\beta=y/m\le 1/401$.
For $y\ge 0.5$, that is $\beta\ge 1/(2m)$. Parametrize by $\beta\in(0,1/401]$ and replace
\begin{align*}
\frac{y+2}{m}=\beta+\frac2m\quad&\text{by}\quad\min(5\beta,\beta+0.001),\\
\frac{2y+3}{m}\quad&\text{by}\quad\min(8\beta,2\beta+0.0015).
\end{align*}
The resulting bound does not depend on $m$.

Over $2\le q\le 400$ and over $y$ or $\beta$ the computed constants yield the following bound for $\Sigma_m(\alpha)/m \leq 0.15784$:
Hence,\[\Sigma_m(\alpha)\le(K-0.0358)m-3\log 2.\]

Combining the two ranges proves the first claim. For the second,
\[
R_m\le e^{\varphi(\theta^*)-0.02208m}\qquad\text{and}\qquad S_m=e^{\varphi(\theta^*)}\sqrt{-2/(\pi \varphi''(\theta^*))}.
\]
Moreover $m^{3/2}e^{-0.02208m}$ decreases for $m\ge 2000$.
\end{proof}

\section{Proof of Theorem~\ref{thm:main}}\label{sec:proof}

It remains to compare consecutive values of $S_m$ and to combine the previous results.

\begin{lemma}
\label{lem:ratio}
For $m\ge 2000$, $S_{m+1}/S_m\ge 1.1$.
\end{lemma}

\begin{proof}
By Proposition~\ref{prop:saddle}(i), $(m+1)\theta_m^*<2\pi$. Since $\theta^*_{m+1}$ maximizes $\varphi_{m+1}$ on $(0,2\pi/(m+1))$,
\[
\varphi_{m+1}(\theta^*_{m+1})\ge\varphi_{m+1}(\theta^*_m)=\varphi_m(\theta^*_m)+\log\bigl|2\sin\bigl((m+1)\theta^*_m/2\bigr)\bigr|.
\]
Here $(m+1)\theta^*_m/2\in\pi[0.790624,\,0.791829\cdot\frac{2001}{2000}]$, where $\sin$ is decreasing, so the last term is at least $\log 2\sin(\pi\cdot 0.791829\cdot 2001/2000)$.

Similarly, by Proposition~\ref{prop:saddle}(ii),
\[
\frac{A_m}{A_{m+1}}\ge\frac{0.129501}{0.131906}\Bigl(\frac{m}{m+1}\Bigr)^3,\qquad\frac{m}{m+1}\ge\frac{2000}{2001}.
\]
Interval arithmetic gives the bound.
\end{proof}

\begin{proof}[Proof of Theorem~\ref{thm:main}]
For $25\le m<2000$ the inequality is verified by exact computation (Section~\ref{sec:comp}). For $m\ge 2000$, Corollary~\ref{cor:twosided}, Proposition~\ref{prop:away} and Lemma~\ref{lem:ratio} give
\[
\frac{M_{m+1}}{M_m}\ge\frac{S_{m+1}}{S_m}\cdot\frac{1-0.005533-3\cdot 10^{-15}}{1+0.004148+3\cdot 10^{-15}}\ge 1.1\cdot 0.9903>1.08.\qedhere
\]
\end{proof}

\section{Computations and formal verification}\label{sec:comp}

The whole proof, including all computations, has been checked in a formalization in Lean~4 with the Mathlib library (version~4.15.0), which is available as supplementary material.

\subsection{Exact part}
For $m\le 2000$ the coefficients are computed in exact integer arithmetic from the recurrence
\[
a_{m,k}=a_{m-1,k}-a_{m-1,k-m},\qquad a_{m,D-k}=(-1)^ma_{m,k},
\]
by a checker whose correctness is proved in Lean. It confirms $M_{m+1}>M_m$ for $25\le m<2000$ and the values $M_0,\dots,M_{26}$, which give $\max_{m\le 25}M_m=16$, $M_{26}=19$ and the inequality used in the proof of Corollary~\ref{cor:conj}.

\subsection{Certified constants}
Every numerical inequality in Sections~2--\ref{sec:proof} is proved in Lean using exact rational arithmetic. Elementary functions are enclosed by Taylor polynomials with proved remainder bounds, and results are rounded outward to multiples of $2^{-64}$ where needed. The finitely many resulting rational inequalities are checked by Lean, the larger ones by compiled code. 

\begin{table}[ht]
\caption{Certified constants (valid for all $m\ge 2000$).}\label{tab:constants}
\begin{tabular}{lll}
\hline
quantity & certified value & reference\\
\hline
$x_0$ & $[0.791224,\ 0.791229]$ & Lemma~\ref{lem:sudler}\\
$K$ & $[0.198607,\ 0.198629]$ & Lemma~\ref{lem:sudler}\\
$\Lambda(\tfrac56)$ & $\le 0.161538$ & Lemma~\ref{lem:sudler}\\
$c^*=m\theta^*/2\pi$ & $[0.790624,\ 0.791829]$ & Proposition~\ref{prop:saddle}\\
$A/m^3$ & $[0.129501,\ 0.131906]$ & Proposition~\ref{prop:saddle}\\
$|\varphi'''(\theta^*)|/m^4$ & $\le 0.057737$ & Proposition~\ref{prop:saddle}\\
$\max_{\mathrm{window}}|\varphi''''|/m^5$ & $\le 0.852654$ & Proposition~\ref{prop:saddle}\\
$M_m/S_m-1$ & $[-0.005533,\ 0.004148]\pm R_m/S_m$ & Corollary~\ref{cor:twosided}\\
gap, $\alpha\le 0.02$ & $\ge 0.02208\,m$ & Proposition~\ref{prop:away}\\
gap, $\alpha\ge 0.02$ & $\ge 0.0358\,m$ & Proposition~\ref{prop:away}\\
$R_m/S_m$ & $\le 3\cdot 10^{-15}$ & Proposition~\ref{prop:away}\\
$S_{m+1}/S_m$ & $\ge 1.1$ & Lemma~\ref{lem:ratio}\\
$M_{m+1}/M_m$ & $\ge 1.0893$ & Section~\ref{sec:proof}\\
\hline
\end{tabular}
\end{table}

\subsection{The Formalization}
The Lean project has about $9{,}500$ lines and follows the structure of the paper, with one file for each of Sections~2--\ref{sec:proof}. Moreover, the TeX file and the Lean file have cross referenced names that one can check. The final theorem of the Lean code is the statement of Conjecture~\ref{conj}. The formalization depends only on the standard axioms of Lean and Mathlib, together with Lean's trust in compiled evaluation for the computational checks. The formal verification has been checked and built successfully using Lean v4.15. The whole codebase is available as a at \href{https://github.com/akuncu/MaxCoeff}{https://github.com/akuncu/MaxCoeff}.

\section*{Acknowledgments}

The authors thank George E. Andrews and Kevin Keating for their interest and for suggesting the circle method as a possible approach toward this problem.

\section*{Declaration of Generative AI Use}
The authors acknowledge the help of Anthropic Claude Opus 5.5 in the computations and drafting of this paper, as well as the formalization provided with the paper. The mathematical content has been driven under the authors' direction, and they take full responsibility for the correctness of the components to the best of their knowledge.

\end{document}